\documentclass[submission, Proceedings]{SciPost}

\hypersetup{
    colorlinks,
    linkcolor={red!50!black},
    citecolor={blue!50!black},
    urlcolor={blue!80!black}
}

\usepackage[bitstream-charter]{mathdesign}
\DeclareSymbolFont{usualmathcal}{OMS}{cmsy}{m}{n}
\DeclareSymbolFontAlphabet{\mathcal}{usualmathcal}

\begin{document}

\begin{center}{\Large \textbf{
Two entangled and scientifically impactful lives:
Ji\v{r}\'i Patera, Pavel Winternitz and the Montréal School of Mathematical Physics
}}\end{center}

\begin{center}
Luc Vinet\textsuperscript{1$\star$}
\end{center}

\begin{center}
{\bf 1} IVADO and Centre de Recherches Math\'{e}matiques, \\
Universit\'e de Montr\'eal, Montr\'eal, QC, Canada\\

* vinet@crm.umontreal.ca
\end{center}

\begin{center}
\today
\end{center}


\definecolor{palegray}{gray}{0.95}
\begin{center}
\colorbox{palegray}{
  \begin{tabular}{rr}
  \begin{minipage}{0.1\textwidth}
    \includegraphics[width=20mm]{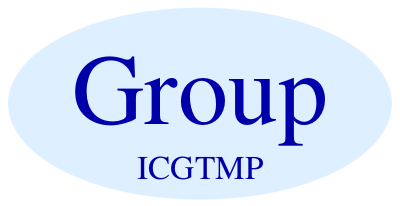}
  \end{minipage}
  &
  \begin{minipage}{0.85\textwidth}
    \begin{center}
    {\it 34th International Colloquium on Group Theoretical Methods in Physics}\\
    {\it Strasbourg, 18-22 July 2022} \\
    \doi{10.21468/SciPostPhysProc.?}\\
    \end{center}
  \end{minipage}
\end{tabular}
}
\end{center}

\section*{Abstract}
{\bf
This text offers a personal account of the scientific legacy of two giants of mathematical physics at the turn of the Millenium and their heritage in Canada, their land of adoption.

\bigskip
   \begin{center}
\textit{To the memory of my mentors and friends, Ji\v{r}\'i and Pavel, with admiration and gratitude.}
   \end{center}
}

\vspace{10pt}
\noindent\rule{\textwidth}{1pt}
\tableofcontents\thispagestyle{fancy}
\noindent\rule{\textwidth}{1pt}
\vspace{10pt}

\section{Introduction}
\label{sec:intro}
Within a year, sadly, Pavel Winternitz and Ji\v{r}\'i Patera passed away in 2021 and 2022. Each one of them stands tall and deserves separate praise. Fate has had it that they have often been celebrated together even though in fact, most of their scientific work has been done separately. I will also indulge in this conflation that will not fully do them justice. One reason is that they cannot be dissociated on the occasion of the 50th anniversary of the International Colloquium on Group Theoretical Methods in Physics (ICGTMP); another, is that they leave a truly joint heritage in Montreal. 

\begin{figure}[h]
	\centering
	\includegraphics[width=0.7\textwidth]{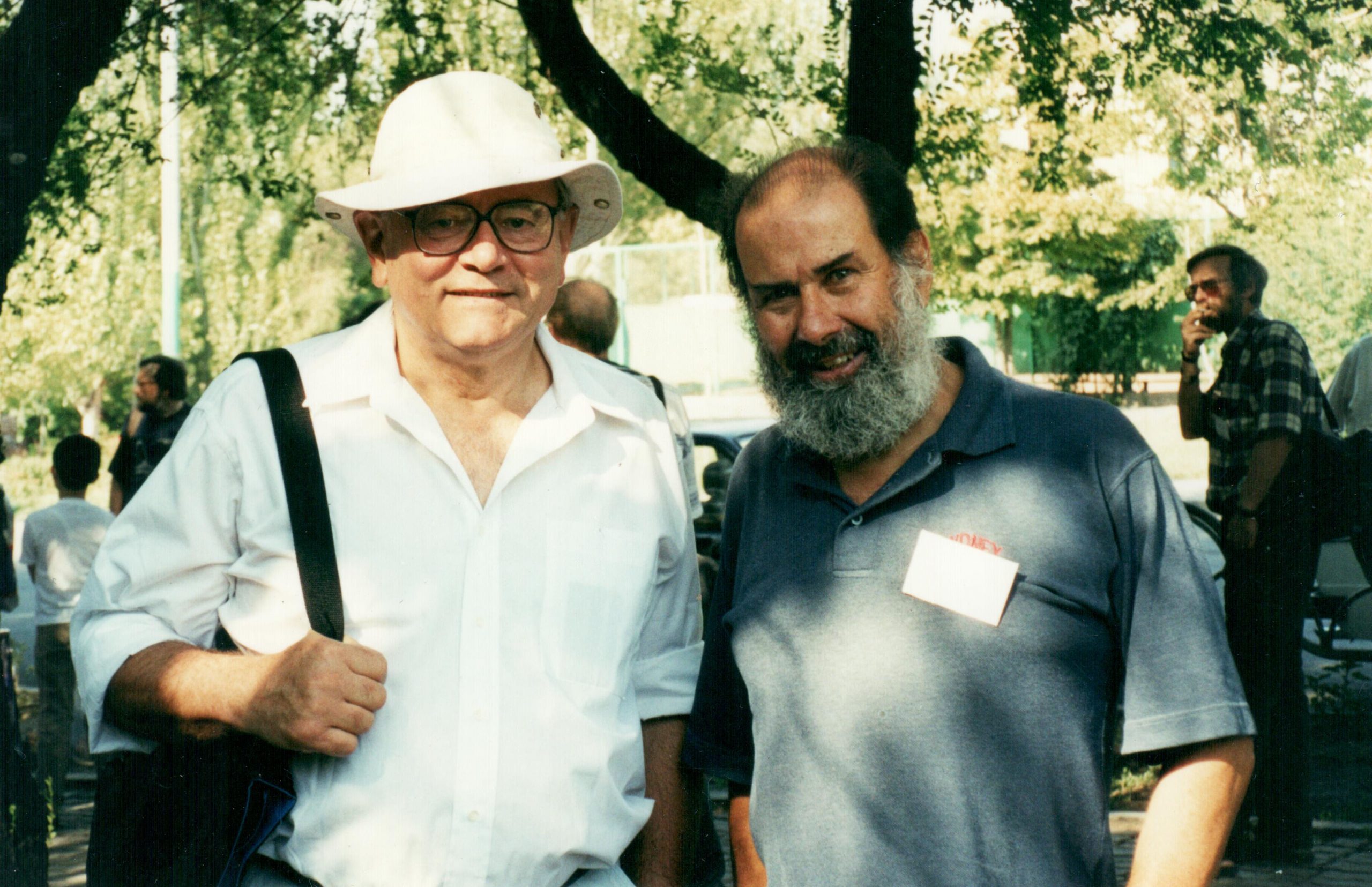}
	\caption{Ji\v{r}\'i and Pavel}
	\label{Centre de recherches mathématiques}
\end{figure}

Like quantum interacting particles that become entangled and yield ``magic" \cite{oliviero2022measuring}, through the vicissitudes of History, these brilliant individuals were set on a colliding course which they steered to have profound influences on various fronts. This is the story that I will try to tell as a tribute to their accomplishments.

\section{Early life and education in the eastern bloc}
Born in Czechoslovakia, both in 1936, Ji\v{r}\'i and Pavel have been educated as theoretical physicists in the great Russian tradition of the Soviet era. 

Ji\v{r}\'i was a native of Zdice, a small town in central Bohemia near Prague. He attended high school in D\v{e}\v{c}\'in in northwestern Bohemia after which he studied at Moscow State University and Dubna as well. In 1964, he obtained his Doctorate from Charles University in Prague. One of his first papers published in 1963 in Nuclear Physics studied the production of $\Lambda$-hyperons in $\pi^- - p$ interactions \cite{barashenkov1963pole} and was written in collaboration with the prominent physicist Blokhintsev, a student of Tamm, who founded the Joint Institute for Nuclear Research (JINR) in Dubna and was its first director. In 1965, Ji\v{r}\'i took a leave from the Physical Institute of the Czechoslovak Academy of Sciences to hold a postdoctoral fellowship from the National Research Council (NRC) of Canada within the developing theory group of the Physics Department of the Universit\'e de Montr\'eal that included Asok Bose, Guy Paquette and soon thereafter Jean-Robert Derome as well as Robert Brunet in the Mathematics Department. This is when, with Bose, Ji\v{r}\'i began his work on group theoretical methods and as you appreciate, this stay in Montreal was to have in many other ways a determining effect on the future. It should be recalled that a memorable world fair the ``Expo 67" took place in Montreal during that period and that the Tchechoslovak pavilion was one of the most popular. At the end of his NRC award in December 1966, Ji\v{r}\'i came back to Prague.

Pavel was born in Prague. He spent the war years in England from where came his fluency in English. He pursued graduate studies at the Leningrad University where Fock was teaching and then in Dubna where he obtained his doctorate in 1966 under the supervision of Smorodinsky, a student of Landau, who Blokhintsev had invited in 1956 to become the head of the Theoretical Group of the JINR. In 1967, he took a leave from the JINR and from the position he had obtained at the Nuclear Research Institute of the Czechoslovak Academy of Sciences in \v{R}e\v{z} to spend time at the International Center for Theoretical Physics in Trieste. In 1968, he returned to his home country.

\section{Prague and Montreal}
Czechoslovakia was an exciting place to be at the beginning of 1968 with the election in January of Dub\v{c}ek as First Secretary. Followed the Prague Spring with its waves of proposed reforms. This unfortunately displeased the Soviet leaders to the extent that in late August troops from four Warsaw Pact countries invaded and controlled Czechoslovakia. For some days, it was possible to flee and so did Ji\v{r}\'i and Pavel. With a visa in hand to attend a scientific meeting in Vienna, Ji\v{r}\'i, his wife Tania and their baby daughter Sacha left with what they could pack in their small car. After a stay in London, they headed to a city that was familiar to Ji\v{r}\'i, namely Montreal, where he was appointed Researcher in the nascent CRM in 1969.
More hesitant to cut ties with Europe, Pavel initially went to England with his wife Milada and their twin boys Michael and Peter and spent some time
at the Rutherford High Energy Laboratory in Chilton where Roger Philips a student of Dirac was his host. A year later he crossed the Rubicon and at the invitation of Wolfenstein, Pavel moved to Pennsylvania with his family first to Carnegie-Mellon and subsequently to the University of Pittsburgh. 
This is when he was encouraged to come to Montreal by Ji\v{r}\'i.

In 1968, the Rector of the Université de Montréal Roger Gaudry, together with Maurice Labb\'e, the first Vice-Rector Research of the university and, Jacques St-Pierre who created the Computer Science Department, had the vision to establish a national institute for research in the mathematical sciences: the Centre de Recherches Math\'ematiques that was to become internationally known as the CRM. They obtained sizable funding from the NRC to that end. Jacques St-Pierre acted as Interim Director and hired Ji\v{r}\'i. Three years later, in 1972, Pavel was also joining the CRM as Researcher. And this is how the stage was set for two young Czechs in their mid-thirties to shape the course of mathematical physics in Canada. In retrospect, the Czech diaspora generated by the 1968 events had a profound effect on the scientific life of Montreal. I shall expand on the roles that Ji\v{r}\'i and Pavel played but there is another striking example that I wish to mention. Montreal has two university hospital systems attached to the Université de Montr\'eal and  McGill University. Quite strikingly at the same time in the 90s and 2000s, the research institutes of both these university hospitals were led by two Czechs: Pavel Hamet and Emil Skamene who had come to Montreal in circumstances similar to those of Ji\v{r}\'i and Pavel and who all became friends of course.

I started my undergraduate studies at the Université de Montr\'eal in 1970. Little was I suspecting that the Prague Spring demise that I had watched unfold with distress two years before was to have a defining impact on my life. In the Fall of 1972, totally oblivious to the creation of the CRM, I wished to enroll in the Master's program in Physics and was looking for a supervisor. It is Robert Brunet from whom I had taken a class who informed me of the existence of the CRM and of the fact that it had recruited two outstanding theoretical physicists with one, Pavel, that had just arrived.  He thought they would be interested in taking graduate students and suggested that I approach them. I followed up and still recall the enthusiasm I felt when together Ji\v{r}\'i and Pavel presented their research programs to me. This was obviously a personal defining moment.

\section {CRM: the early period and the collaborative years}

When the Centre de Recherches Math\'ematiques or CRM was created the plan was to develop research groups in mathematics and statistics, theoretical physics and computer science. While quality hires were made in all three sectors, somehow the group in computer science dispersed. The physics division had a statistical mechanics section that proved less cohesive and very much thanks to Ji\v{r}\'i and Pavel a tradition in mathematical physics developed along with other more mathematical areas that also benefited from the presence of Czech scientists like Anton Kotzig and Ivo Rosenberg.

In the Summer of 1974, when I reported to begin my Master's, the CRM was located in the J\'esus-Marie Pavilion and moved that very Summer to the location on C\^ote Ste-Catherine that many of you have visited and where it stayed until 1994. For starters, I was asked to read the book of Naimark on the representations of the Lorentz group. Eventually, I collaborated more closely with Pavel on superintegrable models. By then, Ji\v{r}\'i and Pavel had hired a postdoc whose name was Ernie Kalnins that I had a hard time understanding because of his New-Zealander accent. The group was already very lively. Bob Sharp who was on Faculty at McGill University had found like-minded colleagues in Ji\v{r}\'i and Pavel at the CRM and was already collaborating with them.  Marcel Perroud who had completed a thesis with Derome and John Harnad became soon regular members of the team whom Ji\v{r}\'i and Pavel have much supported. It is during that period in 1973, that Ji\v{r}\'i together with David Sankoff another CRM researcher, produced, in book format\cite{patera1973tables}, his first set of tables collecting data on representations of simple Lie algebras (others were produced later\cite{mackay1981tables,bremner1985tables,mckay1990tables,kass1990affine}). Thanks to particle physicists such as Pierre Ramond and Dick Slansky who made them largely known, these tables found their ways to the bookshelves of many scientists making practical use of representation theory.

Andr\'e Aisenstadt, a Montreal philanthropist who held a doctorate in theoretical physics from the University of Zurich was a benefactor of the CRM and endowed a distinguished lectureship at the CRM known as the Aisenstadt Chair. In 1973-74, Ji\v{r}\'i and Pavel arranged for Marcos Moshinsky to hold this Chair and to spend an extended period in Montreal. Charles Boyer and Kurt Bernardo Wolf came with him from Mexico on this occasion. Willard Miller Jr. was spending a sabbatical at the CRM (that led to his career-long partnership with Kalnins). You can imagine the intellectual intensity that such a concentration of visiting collaborators was generating and this came to be the rule with Ji\v{r}\'i and Pavel acting as extraordinary magnets. A continuous flow of distinguished speakers would participate in the weekly seminar. Viktor Kac for example who in 1977 had recently completed the classification of Lie superalgebras visited the CRM very soon after he arrived at MIT. In 1979, the Aisenstadt Chair was attributed to Yuval Ne'eman and so on. 
Eugen Dynkin was also among the numerous distinguished people to hold that chair; Ji\v{r}\'i had early on carefully studied his works in Russian and I always thought that this had a big influence on him.

After having completed my Master's, I embarked on a Ph.D. Together with John Harnad and Steven Shnider a differential geometer from McGill and under the benevolent eye of Pavel, I began investigations bearing on gauge field theories and their geometrical and topological properties that were generating much interest at the time. Yvan Saint-Aubin started his Ph.D. two years after I did and other talented students kept coming to this exceptional environment that Ji\v{r}\'i and Pavel were animating.

Between 1973 and 1980, Ji\v{r}\'i and Pavel wrote an astonishing average of 7 papers per year together. The Journal of Mathematical Physics or JMP which was then one of their favorite venues could as well have been called the JJP for journal of Ji\v{r}\'i and Pavel! Their first paper published in 1973 \cite{patera1973new} brought Heun polynomials into the realm of the rotation group representations and 50 years later is still inspiring new results. Ji\v{r}i and Pavel had the knack for developing lasting collaborations with truly distinguished researchers. In the 70s and the 80's one of those was Hans Zassenhaus who is known as a pioneer of computer algebra. He has had an illustrious career that began by being an assistant to Emil Artin in Hamburg in 1936. After having occupied positions at various universities with McGill among those, in 1965 Zassenhaus settled at Ohio State University for the rest of his career and from there visited Montreal regularly even taking a sabbatical at the CRM in 1977-78. As an illustration of the fruitful collaboration he enjoyed with Ji\v{r}i and Pavel I will recall the program they initiated in 1975 aiming to determine the continuous subgroups of the fundamental groups of physics \cite{patera1975continuous}. This involved many additional collaborators (among them Guy Burdet and Martine Perrin from Marseille for work on the optical group), and culminated with the study of the conformal group from that perspective. The results found many applications among which the classification by Beckers (from Li\`ege), Harnad, Perroud and Winternitz of the tensor fields invariant under conformal transformations \cite{beckers1978tensor}. My thesis work on solutions of the Yang-Mills equations has also roots in these foundational studies. And thus Ji\v{r}\'i and Pavel thrived and drew many into their wake.

From 1973 to 1982, the Director of the CRM was Anatole Joffe. At the end of his term, the funding model of the CRM was modified to one where its researchers would hold Faculty positions in university departments and it is then that  Ji\v{r}\'i and Pavel became professors in the Department of Mathematics and Statistics of the Universit\'e de Montr\'eal. In 1984, Francis Clarke was appointed Director and the CRM blazed new trails under his leadership. It played a pioneering role together with the MSRI in defining in the 80s the modern organization of research in the mathematical sciences around institutes organizing visitors and thematic programs.

\section{ICGTMP and outreach}

The International Colloquium on Group Theoretical Methods in Physics was initiated by Henri Bacry from Marseille and Aloysio Jenner from Nijmegen and oscillated between these two cities from 1972 to 1975. In 1974, even though I was a rather young graduate student then, I had the privilege to attend the third edition of the colloquium in Marseille. There was a good contingent from Montreal led of course by Ji\v{r}\'i and Pavel. I could thus observe firsthand the interest they were generating and how many interactions they were having. During the event, they proposed that the conference be held at the CRM in 1976, the year Montreal was to have the Olympic Games. This proposal was adopted and after returning to Nijmegen in 1975, the ICGTMP was held outside Europe for the first time and, thanks to the leadership of Ji\v{r}\'i and Pavel, became truly an international series. Before the colloquium, they also organized a three-week summer school in the framework of the yearly S\'eminaire de Math\'ematiques Sup\'erieures (SMS) initiated by Maurice Labb\'e in 1962 at the Universit\'e de Montr\'eal. Renowned scientists such as Feza G\"{u}rsey, Sigurdur Helgason, Peter Lax, Louis Michel, Willard Miller, Marcos Moshinsky, and many others lectured at these meetings. Ji\v{r}\'i and Pavel showed us the way. Following in their footsteps Yvan Saint-Aubin and I organized the ICGTMP again in Montreal in 1988 together with an edition of the SMS. Because of the glasnost and perestroika taking place in the Soviet Union, it proved possible to host for the first time in the West many distinguished Russian scientists such as Yuri Manin or Sacha Zamolodchikov who were playing key roles in the development of quantum groups and conformal field theories. I shall mention another ICGTMP story. In 1978-79, having mostly finished my Ph.D. work, I accompanied Pavel for a year in Paris as he was spending a sabbatical at Saclay. I shall come back to this later. During that year, Andr\'e Aisenstadt offered me a scholarship that made it possible to attend the Jerusalem Einstein Centennial as well as the ICGTMP that followed in Kyriat Anavim. This offered memorable experiences with Ji\v{r}\'i and Pavel who was kind enough to share a room with me. Pavel was for many years a member of the ICGTMP Standing Committee and in 2018 he received the Wigner medal.

In addition to the ICGTMP, Ji\v{r}\'i and Pavel organized numerous meetings over the years thus influencing research directions broadly and much contributing to the animation of the international scientific community and the visibility of the CRM. One striking example is the workshop on Symmetries and Integrability of Difference Equations that I ideated with Pavel in 1994 and which Decio Levi helped put in place; as you know, this led to the ongoing SIDE series of biennial conferences. Another example involving Ji\v{r}\'i this time is the thematic program entitled Aperiodic Long Range Order that he and Bob Moody organized albeit at the Fields Institute in 1995 \cite{patera1998quasicrystals}.

Ji\v{r}\'i and Pavel have also been very active in developing international collaborations. To that end they made good use of agreements between France and Belgium; this led in particular to the appointment of V\'eronique Hussin at the Universit\'e de Montr\'eal. Over time Ji\v{r}\'i concentrated more in North America. He developed a very fruitful and long lasting collaboration with Bob Moody who was based in Saskatchewan and Alberta. In 1983-84, he spent a sabbatical at Caltech. Around that time, he began collaborating with Gordon Shaw and became involved in the MIND Research Institute which was created in 1998 in Irvine. Ji\v{r}\'i has also been a regular participant in the Aspen Center for Physics program and often visited the MSRI. From the mid 90s onward, with Eliza Shahbazian, Ji\v{r}\'i also pursued collaborative projects with Lockheed Martin Canada and OODA technologies. The international collaborations of Pavel were concentrated mostly in Europe, more precisely in Italy, with Decio Levi and others and in Spain, especially with Miguel Angel Rodriguez and Mariano Del Olmo who had been postdocs at the CRM. He also had close ties with Mexico, with Sacha Turbiner in particular. Furthermore, when it became possible both Ji\v{r}\'i and Pavel reconnected with their roots and put in place collaborative links between Prague and Montreal. This brief overview of the outreach activities of Ji\v{r}\'i and Pavel is grossly incomplete but hopefully illustrates how they have both made Montreal an international hub of mathematical physics.

\section{Moving in different scientific directions: relentless creativity}

As already indicated, Pavel took a sabbatical year in 1978-79 while for Ji\v{r}\'i this happened in 1983-84. At the time, together, they had more or less completed the large undertakings described before and this had monopolized them fully. Without turning their backs on these programs, they then wished to explore new directions and used the occasion of their leaves to do so. As a result, while they kept writing joint papers until 1999 on the classification of maximal Abelian subalgebras \cite{patera1983maximal} and graded contractions \cite{couture1991graded} in particular, the intensity of their collaborative production diminished as each one of them independently opened up new domains. In France, while pursuing his never-ending interest in the nucleon-nucleon scattering phenomenology \cite{bystricky1978formalism}, Pavel decided to focus his attention on the field of non-linear integrable systems whose study with emphasis on solitonic waves and the introduction of the inverse scattering method had been generating great advances. As for Ji\v{r}\'i, his collaboration with Bob Moody had already kicked off with papers on weight multiplicities \cite{moody1982fast} and on characters of elements of finite order\cite{moody1984characters} . This would launch four decades of pioneering research by the two of them of which I will only give a succinct overview \footnote{I shall cite below a number of papers for illustration's sake. This is somewhat arbitrary of course, and I apologize to important collaborators that are not mentioned.}.

From the late 80's, with several collaborators Ji\v{r}\'i developed the profound theory of Lie gradings that he had initiated with Zassenhaus \cite{patera1989lie}. Joris Van der Jeugt who had collaborated with Bob Sharp held a NSERC Visiting Researcher position at the CRM in that period and got involved in those studies. Then, together with Moody, Ji\v{r}\'i made fundamental advances toward the mathematical understanding of quasicrystals viewed as cut and project point sets \cite{moody1993quasicrystals, patera1998quasicrystals}. This led them, while Ji\v{r}\'i was holding a Killam scholarship, to study Voronoi domains \cite{moody1995voronoi} and non-crystallographic root systems \cite{chen1998non}. Ji\v{r}\'i also pursued the applications of these results in cryptography \cite{patera2008quasicrystals}. One additional broad topic that Ji\v{r}\'i has much shaped with signal processing in mind, is that of orbit functions. He wrote a foundational paper with Anatolyi Klimyk \cite{klimyk2006orbit} and much collaborated on this with Ji\v{r}\'i Hrivn\'ak \cite{hrivnak2009discretization} who was a postdoc at the CRM and is now on Faculty at the Czech Technical University in Prague. As part of this program, a generalization of the known properties of the Chebyshev polynomials of the second kind in one variable to polynomials of many variables based on the root lattices of compact simple Lie groups of any type and  any rank was provided \cite{moody2011cubature}. Another fascinating application that Ji\v{r}\'i has explored is the connection that non-crystallographic Coxeter groups have with fullerene and nanotube structures \cite{bodner2013c70}.

Pavel's first contributions to integrable models had to do with B\"{a}cklund transformations. Renewing with Bob Anderson an acquaintance from the time in Trieste, together with John Harnad, he determined the nonlinear superposition properties of matrix Riccati equations \cite{harnad1983superposition}. Subsequently, he launched a broad program aimed at obtaining solutions to these integrable nonlinear partial differential equations through symmetry reduction. This involved finding first the symmetry algebra of the system, a task that he computerized with colleagues \cite{champagne1991computer} and second,  imposing invariance under subalgebras with the help of his expertise at classifying those. This was applied fruitfully to many systems and in particular to the KP one \cite{david1986symmetry} in collaboration with Daniel David, a Ph.D. student of Pavel, Niky Kamran a postdoc at the time and now a distinguished Faculty at McGill and  Decio Levi with whom Pavel wrote the largest number of papers. Pavel and Decio further introduced the notion of conditional symmetry to treat analogously the Boussinesq equation for example \cite{levi1989non}. Michel Grundland who came to the CRM from Poland in the early 80s also participated in these studies. Meanwhile, Pavel pursued his maximal Abelian subalgebras  program (see for example\cite{hussin1990maximal})  as well as the one aimed at characterizing Lie algebras \cite{rand1988identification}, an undertaking that Libor \v{S}nobl joined as a postdoc \cite{vsnobl2005class} to eventually bring it all together in a book \cite{vsnob2017classification} co-authored with Pavel. Another major accomplishment of Pavel has been to develop the Lie theory of difference equations \cite{levi1991continuous, levi1997lie}; the large body of results he has obtained in this area has been collected in a book \cite{levi2022continuous} written with Levi and Yamilov that will be posthumously published. One cannot write about Pavel's scientific production without mentioning his work on superintegrable models which is rooted in his seminal paper \cite{frivs1967symmetry} of 1966. Throughout his career, he kept returning with numerous co-workers to this fertile topic that he championed connecting it for instance to separation of variables and more lately to Painlev\'e transcendents. Pavel's discovery in 2009 with his student Tremblay and Turbiner of the so-called TTW model \cite{tremblay2009infinite} exhibiting constants of motion of arbitrary degrees had the effect of a bomb and  gave an enormous impetus to the field. His former student Ian Marquette and former postdocs Sarah Post and Adrian Escobar-Ruiz among others worked actively with Pavel on this topic in more recent times.

Without adequately summarizing their abundant research outputs, I trust this short overview nevertheless gives a sense of the diversity, richness, and importance of their work.

\section{Conclusion : Passion for research - Legacy and memories}
	
	Ji\v{r}\'i and Pavel both had an unquenchable passion for research and science which they followed with talent throughout their professional lives for more than 60 years. They had the good fortune to never lose their vivacity and their curiosity remained high and sharp. Unfettered by fashions, they pursued their interests to achieve bodies of work of great depth and originality. Even afflicted by blindness in his later years, Ji\v{r}\'i admirably carried on serenely, supervising in this period many students and postdocs who came from the Czech Republic. 
	
	At the Universit\'e de Montr\'eal the title of Emeritus Professor is a high distinction and is awarded parsimoniously; there is a yearly quota for these appointments and nominations across the university are carefully assessed. To be eligible, you need of course to announce your retirement. The title of Emeritus Professor was bestowed upon Pavel in June 2020. That this recognition made him very proud reflects how much his belonging to the CRM and the Universit\'e de Montr\'eal mattered to him. As for Ji\v{r}\'i, he did not find the time to retire!
	
	I had the privilege to be the Director of the CRM from 1993 to 1999  and from 2013 to 2021. This means that I was at the helm in 1996 and 2016 when Ji\v{r}\'i and Pavel turned 60 and 80 respectively. For highly distinguished colleagues to whom much is owed, it is a nice and appropriate tradition to organize celebratory events to express esteem and gratitude. Ji\v{r}\'i and Pavel deserved such an homage and we made sure not to miss out on this occasion in 1996 and so Yvan Saint-Aubin and I organized a conference entitled Algebraic Methods in Physics: A Symposium for the 60th Birthdays of Ji\v{r}\'i Patera and Pavel Winternitz \cite{saint2012algebraic}. For obvious reasons, this was one of those instances where they were lauded together. Present at this event and no longer with us were the Wigner medal recipients Louis Michel, Marcos Moshinsky, and Lochlain O'Raifeartaigh as well as David Rowe and Dick Slansky. Jean-Pierre Gazeau, Basil Grammaticos, Ronald King, Frank Lemire, George Pogosyan, Peter Olver, Alfred Ramani, Guy Rideau, Keti Tennenblat and Ji\v{r}\'i Tolar whom I have not mentioned before were among the participants. 
	
	\begin{figure}[h]
		\centering
		\includegraphics[width=0.7\textwidth]{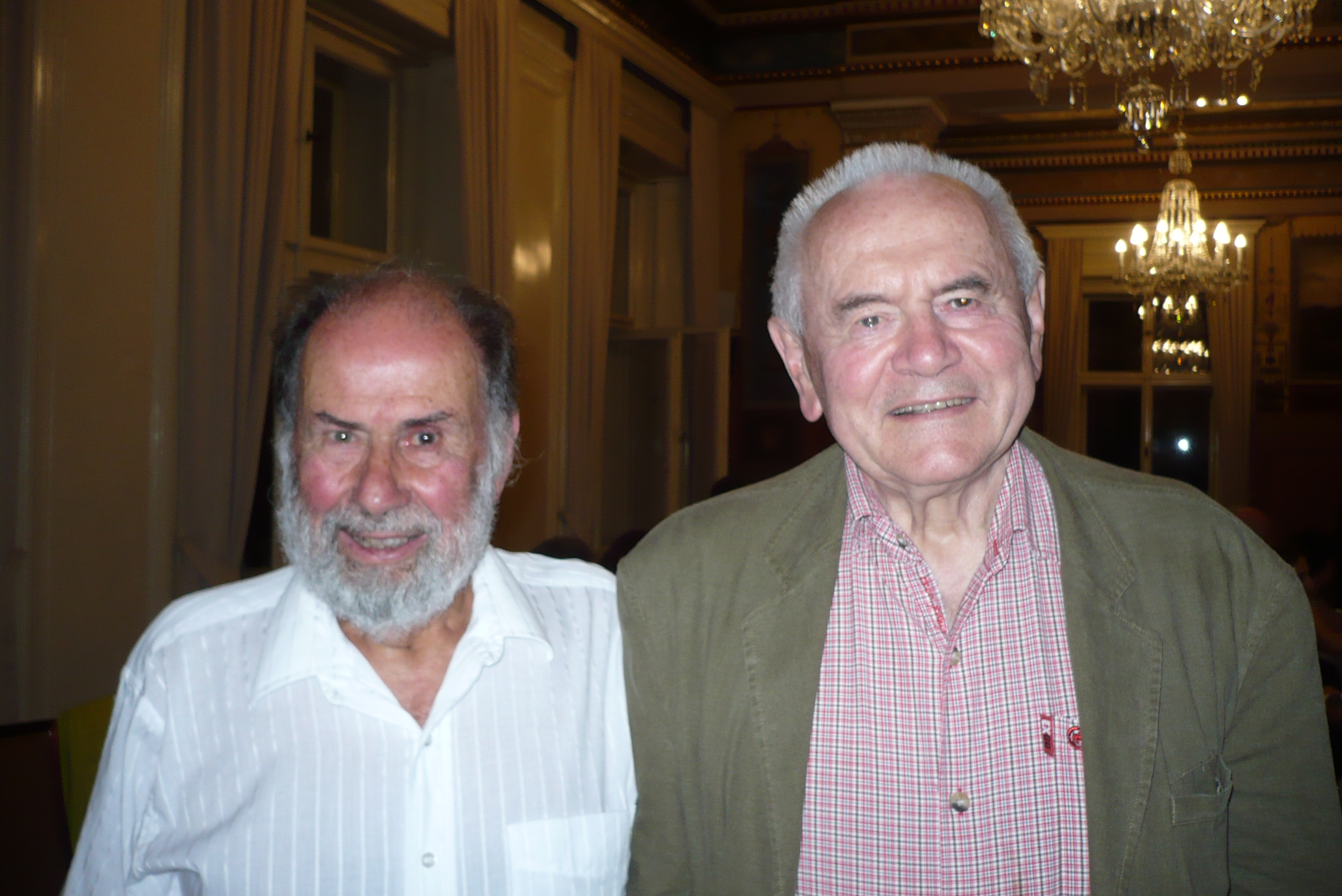}
		\caption{80th Birthdays of Ji\v{r}\'i and Pavel held at Prague in 2016 }
			\end{figure}

	In 2016, both 80 then, Ji\v{r}\'i and Pavel were active as ever and many new generations of colleagues had profited from their interactions with them. Another birthday party was thus in order. This one was held in Prague and put together by the Doppler Institute in collaboration with the CRM. Many Czech colleagues obviously attended, among them Pavel Exner and Igor Jex the Dean of the Faculty of Nuclear and Physical Engineering at the Czech Technical University that was hosting the meeting nicely organized by Libor \v{S}nobl. Let me also add the names of more friends of Ji\v{r}\'i and Pavel who spoke on this occasion and who had not appeared in these lines yet: Vladimir Dorodnitsyn, Hubert De Guise, Luigi Martina, Anatoly Nikitin, Alexei Penskoi, Marzena Szajewska, Piergiulio Tempesta, Mark Walton.\footnote {I have made the perilous choice to identify for memory many individuals whose paths crossed those of Ji\v{r}\'i and Pavel one way or the other. To the many who have unfortunately been left out I apologize trusting they will have understood the intent and will not hold grudges.} It was great to have these occasions to express to Ji\v{r}\'i and Pavel during their lifetime our deep appreciation for their science, their friendship and the bridge they built between Prague and Montreal.

	Alas, they are no longer with us but their legacy lives on. On the scientific front, they have written papers that will keep being touchstones for major areas of mathematics and theoretical physics as well as springboards for many discoveries to come. To all the people they have trained, inspired and befriended, sharing their knowledge, intelligence and culture, they have offered something of themselves that will be transmitted through generations. And, on the human side, they left us with the memory of kind and free men, of proud Czechs and Canadians who were citizens of the World, of fellows who enjoyed life and were caring, of  extraordinarily hospitable and welcoming individuals who instilled in the CRM warmth and excellence and taught it to always aim higher. They have set the stage for the members of their Laboratory to carry on and for other researchers from the world over to join the CRM, walk in their footsteps and like them have a global impact. Ji\v{r}\'i and Pavel, here is to you.

\section*{Acknowledgements}
First, I wish to thank Ji\v{r}\'i Patera and Pavel Winternitz for all they have done for me. I also want to acknowledge the precious help provided by Jos\'ee Savard, Sacha Patera and Michael and Peter Winternitz in the preparation of this presentation.


\paragraph{Funding information}
	My graduate studies in the mathematical physics group of the CRM have been made financially possible thanks to a scholarship from the Natural Sciences, and Engineering Council (NSERC) of Canada. Without this bursary my journey with Ji\v{r}\'i and Pavel might not have happened. Like me, both Ji\v{r}\'i, and Pavel have continuously held NSERC discovery grants. These allow to provide stipends to deserving students or postdocs, and to carry on the training tradition that Ji\v{r}\'i and Pavel initiated. I am sure they would join me in acknowledging this unwavering support. We have also benefited from grants provided by the Fonds de Recherche du Qu\'ebec (FRQ).  Throughout the years, NSERC and FRQ have also kept funding this marvelous institute for research in the mathematical sciences that the CRM is, that Ji\v{r}\'i and Pavel contributed so much to and that I was honoured to lead for 14 years. In the name of Ji\v{r}\'i and Pavel, I wish to thank all the Canadian taxpayers who, through the various funding agencies, are making scientific treks like theirs possible with all the benefits that these generate.






\bibliographystyle{SciPost_bibstyle} 
\bibliography{ref_PW.bib}

\nolinenumbers

\end{document}